\documentclass[twoside,leqno,fleqn]{article}
\usepackage{amsmath}
\usepackage{amssymb}
\usepackage[T1]{fontenc}
\usepackage[shortlabels]{enumitem}
\usepackage[onehalfspacing]{setspace}
\usepackage{calligra}

\newcommand{\AUTOR}{C. Ga\ss ner}
\newcommand{\TITEL}{Permutation models of second order}
\def\S{\Sigma}
\def\I{{\cal I}}
\def\G{\mathfrak{G}}
\def\H{\mathfrak{H}}

\def\F{\mbox{\scriptsize\calligra T\!\!\!\!\!\hspace{0.005cm}F\!\!\!\!\!T\!\!\!\!\!\hspace{0.02cm}F}\,\,\,}

\newcommand{\bbbn}{\mathbb{N}} 
 
\newcommand{\bbbq}{\mathbb{Q}} 

\newcommand{\mbm}[1]{\mbox{\boldmath{$#1$}}}
\newcommand{\mbmss}[1]{\mbox{\scriptsize\boldmath{$#1$}}}
\newcommand{\mbmty}[1]{\mbox{\tiny\boldmath{$#1$}}}
\newcommand{\qed}{\hfill{$\Box$}} 

\newtheorem{satz}{Satz}[section] 
\newtheorem{lemma}[satz]{Lemma}
\newtheorem{proposition}[satz]{Proposition}

\newtheorem{theorem}[satz]{Theorem}
\newtheorem{example}[satz]{Example}
\newtheorem{defi}[satz]{Definition}

\newtheorem{remark}[satz]{Remark}
\newtheorem{overview}[satz]{Overview}
\begin{document}
\newcounter{li}

\thispagestyle{empty}
\begin{center} {\Large\bf Permutation Models of Second Order}
\vspace{0.4cm}\\{\bf Christine Ga\ss ner}
 { \vspace{0.1cm}\\Universit\"at Greifswald, Germany, 2024\\ gassnerc@uni-greifswald.de}\\\end{center} 

\begin{abstract} 
 G\"unter Asser (1981) introduced second-order permutation models. In this way, the Fraenkel-Mostowski-Specker method for defining models of ZFA was transferred to a new application area.  To investigate the strength of second-order principles of choice in second-order predicate logic (PLII) with Henkin interpretation (HPL), we have extended the Fraenkel-Mostowski-Specker-Asser method. Here we discuss many details, including some useful references. Finally, we address a question raised by Stephen Mackereth.
\end{abstract}

\setcounter{satz}{0}

\section{Introduction}

The following\, {\em Fraenkel-Mostowski-Specker-Asser   method}\, for constructing permutation models of HPL was  introduced in \cite[pp.\,\,112--114]{1c} and applied  in \cite{Gass84} and \cite{Gass94}.  The concept goes back to  observations that were  used  by  Fraenkel (cf.\,\,\cite{FA22a}) and Mostowski to construct the so-called permutation models  of  ZFA. Various Henkin-Asser structures  can be constructed, for instance, by using  any group  $\G$ of permutations of a non-empty set $I$ of so-called individuals and   any normal filter --- a special system --- of subgroups of $\G$. Each resulting  permutation model is a model of second order consisting of predicates $\alpha$ on $I$ for which $\widetilde \alpha$ is symmetric  in the sense of Fraenkel \cite{FA22a}. The resulting models of HPL have  two  important properties. On the one hand,  the  models are  closed with respect to the permutations in the used group $\G$. This means that the interchange of  individuals via permutations  is possible and  causes an allowed   transformation of   predicates into predicates.  More precisely,   it causes   transformations  of predicates  by permutations in $\G$  in such a way that the  resulting   predicates belong to such a  model iff  the original predicates also belong to the model. On the other hand,  we assume that there exists  a suitable filter so that  a predicate is in such a model iff  there is  a certain subgroup in this filter with the property  that the predicate is invariant   with respect to  all permutations belonging to  this subgroup. By definition, we will consider only normal filters.   This means that a predicate can be in such a model  only if it is self-conjugate in the sense of Fraenkel \cite{FA22a}. Moreover,   each  second-order formula   $H$ describing a relationship between a finite number of  individuals (or atoms) and   predicates (or relations) stays $true$    in such a permutation model  even if some  individuals or  predicates assigned to the free variables  in  $H$ are changed by  a  permutation belonging to a suitable subgroup of the filter.
Consequently, the restriction to normal filters   implies --- in any case --- that second-order permutation models are closed with respect to second-order definability.
 Fraenkel \cite{FA22b, FA22a} constructed  permutation  models of set theory with atoms  to  show the independence of the Axiom of Choice from the axioms of ZFA. One of the constructed models is a permutation model of  ZFA that Thomas J.\,\,Jech  called the {\em second Fraenkel model} in \cite[p.\,48]{Jech73}  (cf.\,\,also Example 15.50 in \cite[p.\,252]{Jech03}). 
The following definitions are closely related to the theoretical foundations of the Fraenkel-Mostowski-Specker  method  as outlined in \cite{Felgner, Jech73, Jech03}. 
In defining models of second order,  we assume that the  theory  ZF  is consistent and we use it as our metatheory. Consequently, we have also the consistency of  ZFA including the assumption that the set of atoms is infinite  (cf.\,\,\cite[p.\,45]{Jech73}).   However, in contrast to the  permutation models of  ZFA where the transitivity of the models is an important property, we do not need  the existence of atoms (in our metatheory).
We can  consider any predicate structure  $(I_n)_{n\geq 0}$  provided that the set $I_0$,  the  relation $\widetilde \alpha$  for  each predicate $\alpha$ in  $I_n$ ($n\geq 1$), the sets $\widetilde I_n$  containing all $\widetilde \alpha$ for all $\alpha \in I_n$ and all  $n\geq 1$,  and  the family $(\widetilde I_n)_{n\geq 1}$  belong to  a model of ZF and ZFC, respectively. 

The preparation of this introduction  to second-order models began with the need to provide a script for the lectures in {\sf Mathematical Logic} in Greifswald about 10 years ago. The second-order logic and second-order models are also considered, i.e.,  in \cite{Henk,1c,Gass84,13, Gass94, Kame,Vaeae,Mackereth}. Here, we will particularly  use \cite{Gass24A}   and the following proposition. Note, that some of the important definitions and propositions are recalled. Here, we consider Henkin-Asser structures. As explained at the CL 2022 (cf.\,\,\cite{Gass22B}) we  distinguish between Henkin-Asser structures, Henkin-Shapiro structures, and Henkin-V\"a\"an\"anen structures.  
\begin{proposition}\label{MeaningOfFreeVar} Let $\S$ be a  structure    in ${\sf struc}_{\rm pred}^{({\rm m})}$, let $H$ be a second-order formula in ${\cal L}^{(2)}$, and let  $x_{i_1}, \ldots, x_{i_n}$ and $A_{j_1}^{n_1},\ldots,A_{j_r}^{n_r}$ be the only   variables that may occur free in $H$. Then,   for all assignments $f$ and $f'$ in  ${\rm assgn}(\S)$ with  $f(x_{i_1})=f'(x_{i_1}),\ldots,  f(x_{i_n})=f'(x_{i_n})$ and  $f(A_{j_1}^{n_1})= f'(A_{j_1}^{n_1}), \ldots,  f(A_{j_r}^{n_r})= f'(A_{j_r}^{n_r})$,  we have 
\[\S_f(H)=\S_{f'}(H).\] 
\end{proposition}

\section{Permutations and automorphisms}

\subsection{Basic structures}

${\sf struc}^{(1)}(I)$  is the class of all one-sorted first-order structures given in the form   $(I;(\mbm{c}_k)_{k\in K} ; (\mbm{f}_l)_{ l\in L}; (\mbm{r}_m)_{m\in M})$) and defined as usual.
An {\em automorphism} of a structure ${\sf S}$  in ${\sf struc}^{(1)}(I)$    is an isomorphism of  ${\sf S}$ onto ${\sf S}$. Let us mention that the only automorphism of a structure containing a well-ordering of $I$  as relation   is the identity permutation  (cf.\,\,\cite[Corollary 2.5, p.\,18]{Jech03}). For any non-empty set $I$, each {\em permutation of $I $} is a bijection from $I$ to $I$. Let $({\rm perm}(I),\circ)$ and $(\G_{\sf 1}^I,\circ)$ (shortly denoted by ${\rm perm}(I)$ or $\G_{\sf 1}^I$) be the group of all permutations of $I$ with the usual function composition $\circ$ as binary group operation. Then, any automorphism of a structure ${\sf S}$ in ${\sf struc}^{(1)}(I)$ belongs to $\G_{\sf 1}^I$. Moreover, all automorphisms of ${\sf S}$ together form a group $({\rm auto}({\sf S}), \circ)$ shortly denoted by ${\rm auto}({\sf S})$. Thus, $ {\rm auto}({\sf S})$ is a subgroup of $\G_{\sf 1} ^I$. For each non-empty set $I$ and the structure $(I;\emptyset;\emptyset;\emptyset)$ that contains neither a constant nor a relation nor a function,  ${\rm auto}(I;\emptyset;\emptyset;\emptyset)=\G_{\sf 1}^I $ holds.

For many-sorted structures ${\sf S}_1$  and ${\sf S}_2$  of signature $\mbm{\sigma}^{(2)}$ given by ${\sf S}_1=(\bigcup_{n\geq 0}I_n;\emptyset;\emptyset; (\mbm{s}_n)_{n\geq 1})$, ${\sf S}_2=(\bigcup_{n\geq 0}J_n;\emptyset;\emptyset; (\mbm{t}_n)_{n\geq 1})$, let  $({\it \Phi}_n)_{n\geq 0} $ be a sequence of bijections ${\it \Phi}_n: I_n \to J_n$. $({\it \Phi}_n)_{n\geq 0} $ is an {\em isomorphism from ${\sf S}_1$ to ${\sf S}_2 $}  if, for every $n\geq 1$, ${\it \Phi}_n$ satisfies  $({\it \Phi}_n(\alpha),{\it \Phi}_0(\xi_1),\ldots,{\it \Phi}_0(\xi_{n})) \in \mbm{t}_n$ for any $\alpha\in I_n$ and $(\xi_1,\ldots,\xi_n)\in I^n$ iff $(\alpha, \xi_1,\ldots,\xi_n)\in \mbm{s}_n$ holds. 
For more details see \cite[pp.\,7--8]{1c}. 
An {\em automorphism} of a structure ${\sf S}$ of signature $\mbm{\sigma}^{(2)}$ is an isomorphism $({\it \Phi}_n)_{n\geq 0} $ of ${\sf S}$ onto ${\sf S}$. All automorphisms of ${\sf S}$ together form a group denoted by ${\rm auto}({\sf S}) $. 

\subsection{The idea}

We will  apply  the following method allowing to construct Henkin structures  with the help of suitable groups of permutations in analogy to the construction  of models of ZFA by the Fraenkel-Mostowski-Specker  method. For constructing  a Henkin structure we will   use a  group of automorphisms of a suitable basic  structure with a domain $I$ of individuals. More precisely speaking, each of the {\em basic structures}, ${\sf S}$, can be an arbitrary  one-sorted first-order structure in $ {\sf struc}^{({\rm 1})}(I)$    or a one-sorted first-order relational  structure $\widetilde\S\in {\sf struc}_{\rm rela}^{({\rm 1})}(I)$ (without functions) derived from a predicate structure $\S$ in ${\sf struc}_{\rm pred}^{({\rm m})}(I)$. For an example see Example \ref{ModelMost}.  Note that in general the  individual domain  $I$ is infinite and, thus, the number  of  all permutations  of this domain is  also   infinite. In every case, we are interested in  predicates on $I$ and the  subgroups  of automorphisms that transform  these predicates into themselves. Given a group $\G$ of automorphisms of ${\sf S}$, we will call a subgroup of  automorphisms  in $\G$ under whose application a predicate on $I$ stays invariant --- as usual --- the  symmetry subgroup of $\G$ with respect to  this predicate (see Section \ref{SectionSymmetrySubgroups}) or the symmetry  subgroup of this predicate (if $\G$ is fixed). We will consider specific systems of subgroups of automorphisms so that  all predicates, each of which has a symmetry   subgroup belonging to such a system,  form  a Henkin structure. Each   normal filter on $\G$ that  is a suitable system of subgroups  for this purpose  can be defined explicitly or  via  a normal ideal on the universe of ${\sf S}$ that allows to define subgroups for generating a normal filter (see Definition \ref{DefNormFilt} and Definition \ref{DefExtNormIdeal}).  Hence, in the following we define Henkin structures  $\S(I, \G,\F)$ by taking a group  $\G$ of  automorphisms of a basic structure  ${\sf S}$ with the universe $I$, using   a normal filter  $\F$ of subgroups of $\G$ (see Definition \ref{DefiPermModelByFilter}) or a normal  ideal for generating $\F$, and collecting all predicates whose  symmetry subgroups are in the corresponding  filter $\F$.   Let us consider more details for ${\sf S}=(I;\emptyset;\emptyset;\emptyset)$ and other simple basic structures.

\subsection{Permutations  and their first extension on predicates}  
Let $I$ be any arbitrary non-empty set of individuals and $\pi$ be a permutation in $\G_{\sf 1}^I$.  Then, for all $\mbm{\xi} \in I^n$, we will use the abbreviation $\pi(\mbm{\xi})$ defined by  $\pi(\mbm{\xi})= (\pi(\xi_{1}), \ldots,\pi(\xi_{n}))$.  Moreover, for each $n$-ary predicate $\alpha\in {\rm pred}_n(I)$, let $\alpha^\pi$ be the  $n$-ary predicate satisfying  \[\alpha^\pi(\mbm{\xi})=\alpha(\pi^{-1}(\mbm{\xi}))\]  for  all $\mbm{\xi}\in I^n$.\footnote{The definition of  $\alpha^\pi$ given in \cite[p.\,112]{1c}  was modified.}   Thus, we have $\widetilde{\alpha^{\pi}}= \{\pi(\mbm{\xi})\mid \mbm{\xi}\in \widetilde\alpha\}$. Let $\G$ be any subgroup  of  $\G_{\sf 1} ^I$. The  group  $(\G,\circ)$ (shortly,  $\G$) can be  any  group of  bijections from $I$ to $I$  or a group  of automorphisms of a given basic structure ${\sf S}$ with $I$ as universe.  It is easy to   extend every  permutation $\pi$  in $\G$ uniquely  to a permutation $\pi^*$ of $I\cup {\rm pred}(I)$ by defining
\[\pi^*(\alpha)=\alpha^\pi \quad \mbox{ for all }\alpha \in {\rm pred}(I)  \quad \mbox{ and }\quad \pi^*(\xi)=\pi(\xi) \quad \mbox{ for all }\xi\in I. \]    
Let   $\G^*= \{ \pi^* \mid \pi\in \G\}$. Consequently,    $\G^*$  is a subgroup of  $(\G_{\sf 1} ^I)^*$.  $\G^*$  and   $(\G_{\sf 1} ^I)^*$ are groups   that act on $I\cup {\rm pred}(I)$ from the left.  Let $\S$ be any predicate structure $(J_{n})_{n\geq 0}$ in ${\sf struc}_{\rm pred}^{({\rm m})}(I)$ and let 
\[\begin{array}{ll}\G^*|_{\S}=&\{   g:\bigcup _{n\geq 0}J_n\to I\cup {\rm pred}(I)  \mid\\&(\exists \pi \in \G)((\forall \xi  \in J_0) (g(\xi)=\pi(\xi))\,\, \&\,\,\,(\forall \alpha \in\bigcup _{n\geq 1}J_n) (g(\alpha)=\alpha^{\pi}))\}\end{array}\]

\noindent be the set of all functions being the  restrictions of the permutations in $\G^*$ to  $\bigcup_{n\geq 0}J_n$.  For  any $\pi\in \G ^I_{\sf 1}$ and any subset $J\subseteq \bigcup _{n\geq 0}J_n$, let $\pi^*|_{J}$   be the restriction of $\pi^*\in (\G ^I_{\sf 1})^*$ to  $J$ and let $ \pi^ *|_\S=\pi ^*|_{\bigcup _{n\geq 0}J_n}$.  By definition, we thus have $ \pi^ *|_\S=\bigcup _{n\geq 0} \pi^ *|_{J_n}$.

Consequently, for $\S=(J_{n})_{n\geq 0}$, a  permutation  $\pi$ of  the domain of individuals  and its extension $ \pi^ *|_\S\in \G^*|_{\S}$ preserve the structural properties of $\S$ and guarantees the closeness of $\S$ with resect to  $\pi$ only if $ \pi^ *|_\S$ is  in $ {\rm perm}(\bigcup _{n\geq 0}J_n)$, which means, only if $( \pi^ *|_{J_n})_{n\geq 0}$ is in  ${\rm auto}(\S)$. Note, that $\pi\in {\rm auto}(\widetilde\S) $  is not necessary. 

Let us define $ {\rm auto}_\G^*(\S)$ and $ {\rm auto}^*(\S)$ by $ {\rm auto}_\G^*(\S)= \G^*|_{\S}\cap {\rm perm}(\bigcup_{n\geq 0}J_{n})$ and $ {\rm auto}^*(\S)= {\rm auto}_{\G ^I_{\sf 1}} ^*(\S)$ for $\S\in {\sf struc}_{\rm pred}^{({\rm m})}(I)$. 

\begin{lemma}[Closeness of $\S$ with respect to $ {\rm auto}_\G^*(\S)$]  Let $\S$ be a predicate structure $(J_{n})_{n\geq 0}$ and $\G$ be a subgroup of $\G_{\sf 1}^{ I}$.  Then,  $ {\rm auto}_\G^*(\S)$ is a group  that  acts on  $\bigcup _{n\geq 0}J_n$.
\end{lemma}

\noindent
{\bf Closeness with respect to permutations.}  Let  $I$ be any non-empty set, $\S$ be any predicate structure  $(J_n)_{n\geq 0}$ in ${\sf struc}^{({\rm m})}_{\rm pred}( I)$, and $\pi\in \G_{\sf 1}^{ I}$. Consequently, $J_0=I$.  We say that the  predicate domain  $J_n$ {\em is closed with respect to  $\pi$} if $\alpha ^\pi $ belongs to $ J_n$ for each  $\alpha\in J_n$. This means that $ \pi^ *|_{J_n}\in {\rm perm}(J_n)$.  The  predicate structure $\S$  {\em is closed with respect to a permutation $\pi\in \G_{\sf 1}^{ I}$} if all its predicate domains $J_n$ ($n\geq 1$) are  closed with respect to $\pi$.  Thus, the permutations with respect to which $\S$ is closed form a semigroup  of permutations in $\G_{\sf 1}^{ I}$.  This means that, for the  usual function composition operation $\circ$ and all $\pi_1,\pi_2\in \G_{\sf 1}^{ I}$, $\S$ is closed with respect to  $\pi_1\circ\pi_2$ if  $\S$ is closed with respect to  $\pi_1$ and to $\pi_2$.  $\S$ is  closed with respect to all permutations $\pi$ in a subgroup $\G$ of $\G_{\sf 1}^{ I}$ only if it is also closed  with respect to their inverses $\pi^{-1}$.

\begin{lemma}[Closeness of $\S$ with respect to $\pi$ and $\pi^{-1}$]\label{Clos_Struc} \hfill
 Let $\S$ be any \linebreak predicate structure $(J_{n})_{n\geq 0}$, $\G$ be a subgroup of $\G_{\sf 1}^{ I}$, and  let  $\pi$ be a permutation in $\G$ such that  $\pi^*|_\S\in {\rm auto}_\G^*(\S)$ holds.  Then, $\S$ is closed with respect to $\pi$ and $\pi^{-1}$.
\end{lemma}

 \vspace{0.1cm}\noindent {\bf The influence of permutations on the truth values of formulas.} 
A useful property of the permutation models of ZFA  is the closeness of every model with respect to the application of an extended permutation (cf.\,\,\cite[pp.\,45--47]{Jech73}). This property is  a useful prerequisite here, too. However, it is important to know that the values that are relevant in defining a predicate  could be changed by the application of an extended permutation. Our  permutation models should be models of  all axioms in $^hax^{(2)}$ and in particular of the axioms given by (A4).  This implies the requirement for closeness of the permutation models with respect to the definability of predicates.   Thus, we have to  take into account how the truth values of second-order formulas are influenced by changing  an assignment in a predicate structure. For this purpose, we extend  the domain  of each permutation in  $(\G_{\sf 1}^I)^*_\S$   to a domain including all assignments in ${\rm assgn}(\S)$. 

\subsection{Permutations  and a second extension to assignment functions} Let  $I$ be any non-empty set, $\S$ be any predicate structure  $(J_n)_{n\geq 0}$ in ${\sf struc}^{({\rm m})}_{\rm pred}( I)$ which means $I=J_0$, and $\pi\in \G_{\sf 1}^I$.  Then, it 
shall be possible to apply  a function   from ${\rm assgn}(\S)$ to  ${\rm assgn}(\S_{\sf 1}^{I})$ determined by $\pi$. For every  $f$  in    ${\rm assgn}(\S)$,   let $f^{\pi}$ be the assignment in ${\rm assgn}(\S_{\sf 1}^{I})$ defined  as follows.\footnote{The definition  of $f^{\pi}$ given in \cite[p.\,112]{1c}  was modified.}   

 \vspace{0.1cm}
\begin{tabular}{ll}

 $f_0^{\pi}(x_i)$ & $=_{\rm df}\pi(f_0(x_i)) $\hfill{ \hspace{2.8cm} (for all  $f_0:\{x_i\mid i\geq 0\}\to I$)}\\[1.5ex]

 $f_n^{\pi}(A_i^n)$ & $=_{\rm df}(f_n(A_i^n))^{\pi} $\hfill{  (for all  $f_n:\{A^n_i\mid i\geq 0\}\to  J_n$, $n\geq 1$)}\\[1.5ex]

 $f^{\pi}$ & $=_{\rm df}(f_n^{\pi})_{n\geq 0}$\hfill{  (for all  $f=(f_n)_{n\geq 0}$)}
\end{tabular}
 \vspace{0.1cm}

\noindent  For  $f_0:\{x_i\mid i\geq 0\}\to I$ we get   $f_0^\pi:\{x_i\mid i\geq 0\}\to I$. For all $f_n:\{A^n_i\mid i\geq 0\}\to  J_n$, we get $f_n^\pi:\{A^n_i\mid i\geq 0\}\to  {\rm pred}_n(I)$.

\begin{lemma}[Closeness of $ {\rm assgn}(\S)$ with respect to $\pi$] Let $\S=(J_{n})_{n\geq 0}$ and $\pi\in \G_{\sf 1}^{J_0} $. If $ \pi^ *|_\S\in {\rm auto}^*(\S)$,  then there holds $f^{\pi}\in {\rm assgn}(\S)$.
\end{lemma}

\noindent  Such an assignment  $f^{\pi}$ has the following three properties   (cf.\,\,\cite[Hilfssatz 1, and Hilfssatz 2 p.\,112]{1c}). 

\begin{lemma}[Asser \normalfont{\cite{1c}}]\label{HSAsser} Let  $\S$  be the predicate structure $(J_n)_{n\geq 0}$ and let  $\pi$ be a permutation in $\G_{\sf 1}^{J_0} $ such that  $ \pi^ *|_\S\in {\rm auto}^*(\S)$ holds. Let $f$ be  an assignment in $\S$ and $H$  be any formula in ${\cal L}^{(2)}$.  Then,  we have {\em (1)}  and   we have the properties {\em (2)} and {\em (3)} if $H$ is in  ${\cal L}^{(2)}_{x_{1},\ldots,x_{n}}$.
\newcounter{liB}
\begin{list}{\normalfont (\arabic{liB})}{\usecounter{liB}\labelwidth0.5cm \leftmargin0cm \itemsep2pt plus1pt
\topsep1pt plus1pt minus1pt
\labelsep4pt \parsep0.5pt plus0.1pt minus0.1pt \itemindent0.8cm}
\item\label{HSAsser1}  $\S_{f^{\pi}}(H)=\S_f(H)$.
\item\label{HSAsser2} $\alpha_{\S,H,\mbmss{x}, f^{\pi}}=\alpha_{\S,H,\mbmss{x}, f} ^{\pi}$.
\item\label{HSAsser3}  $\alpha_{\S,H,\mbmss{x}, f}^{\pi}$ is definable over $\S$  if $\alpha_{\S,H,\mbmss{x}, f} $   is definable over $\S$.  
\end{list}
\end{lemma} 

It is sufficient to prove (1)  by induction on the formula structure and by induction on the depth of PL\,II-formulas, respectively.    
By  property (\ref{HSAsser1}),  we obtain (\ref{HSAsser2}). Thus, we get (\ref{HSAsser3})  by (\ref{HSAsser2}).  

 \vspace{0.1cm}

 By Lemma \ref{HSAsser} (\ref{HSAsser3}), each extension $\pi^*$ of a permutation $\pi$ in a subgroup $\G$ of $\G_{\sf 1} ^I$  transforms  the predicate $\alpha_{\S,H,\mbmss{x}, f}$  definable by $H$,  individuals $\xi_{n+1},\ldots,\xi_{n+s}$, and  predicates $\alpha_{i_1}, \ldots,\alpha_{i_r}$ into  a predicate definable by $H$. If  the mentioned  individuals and predicates are  mapped to  themselves by means of $\pi$, then the relation $\widetilde \alpha_{\S,H,\mbmss{x}, f}$ is also transformed into itself by $\pi$.
Since   we are especially interested in constructing structures satisfying  the axioms given by  (A4),  we consider so-called symmetry subgroups of a group  $\G$ of permutations,   stabilizer subgroups, and  supports for  predicates and discuss their  definability  by means of second-order formulas.

\subsection{Symmetry subgroups and invariance}\label{SectionSymmetrySubgroups}
Let $I$ be any non-empty set and $\G$ be a subgroup of $\G_{\sf 1} ^I$. For each predicate $\alpha \in  {\rm pred}(I)$, let ${\rm sym}_{\G}(\alpha)$ be defined by
\[{\rm sym}_{\G}(\alpha)=\{ \pi \in {\G}\mid  \alpha^{\pi}=\alpha\}\]
  and   let ${\rm sym}_{\G}^*(\alpha)=_{\rm df} ({\rm sym}_{\G}(\alpha))^*$.    Then, ${\rm sym}_{\G}(\alpha)$  is a group   called the {\em symmetry subgroup of  $\G$ with respect  to $\alpha$}. Whereas $\G^*$   is  a subgroup  of $(\G_{\sf 1} ^I)^*$,     ${\rm sym}_{\G}^*(\alpha)$ is  a subgroup  of $\G^*$ that also acts on $I\cup {\rm pred}(I)$ from the left.  $ {\rm sym}_{\G}^*(\alpha)$  is  the  stabilizer subgroup (i.\,e.\,\,the isotropy subgroup) of   $\G^*$ with respect to $\alpha$.   Moreover, let $perm^*(I)=(\G_{\sf 1}^I)^*$.   Thus, for all $\alpha \in {\rm pred}(I)$, we have  $\pi^*(\alpha)= \alpha ^\pi$    for $\pi^*\in perm^*(I)$ and  $\pi\in \G_{\sf 1}^I$, respectively.
In the following, we also omit the asterisk after $\pi$  in $\pi^*$ and write then also $\pi$ if we consider $\pi^*$.

For any subgroup  $\G$ of  $\G_{\sf 1} ^I$ and   any $P\subseteq I\cup {\rm pred}(I)$,  let ${\G}( P)$ be the subgroup  containing  all $\pi \in {\G}$  whose extensions map each element  in $P$  to itself.  Another notation for such a group could --- in accordance with the notations used     in \cite[p.\,47]{Jech73} --- be ${\rm fix}_{\G}(P)$.   More precisely,   let
\[{\G}( P)=\{\pi \in \G\mid \pi(x)=x \mbox{ for all  } x  \in P\}  \mbox{ and  } {\rm fix}_{\G}(P)={\G}( P)\]
 if $P\subseteq I$.  For $P\subseteq{\rm pred}(I)$, let 
\[{\G}(P) = \{ \pi \in {\G}\mid \alpha^\pi= \alpha  \mbox{ for all  }\alpha \in P\}.\] 
Thus, for $P\subseteq I\cup {\rm pred}(I)$, it holds
\[{\G}( P)={\rm fix}_{\G}(P\cap I) \,\cap\,  \G (P\cap {\rm pred}(I)).\] 
 This means that the individuals in $P$ are  fixed points of each $\pi$ in ${\G}( P)$.  Moreover,  ${\G}( P)\subseteq {\rm sym}_{\G}(\alpha)$ is satisfied for each $\alpha \in P\cap {\rm pred}(I)$.  Hence, ${\G}( P)$ is the group of all automorphisms of $(I;P\cap I;\emptyset;(\widetilde\alpha)_{\alpha\in P\cap {\rm pred}(I)})$. If $P\subseteq I\cup {\rm pred}(I)$ is a  set  satisfying  ${\G}( P)\subseteq {\rm sym}_{\G}(\alpha)$ for some $\alpha \in  {\rm pred}(I)$, then we say   that $P$ is a {\em  support of    $\alpha$ with respect to $\G$}.  
If $P\subseteq I$ is a finite set of individuals with the property ${\G}( P)\subseteq {\rm sym}_{\G}(\alpha)$ for some $\alpha \in  {\rm pred}(I)$, then we say  --- by analogy with \cite[p.\,47]{Jech73} --- that $P$ is a {\em finite individual support of   $\alpha$}.  Note, that $\emptyset$ is a support for $\alpha$ with respect $\G(\{\alpha\})$. With the help of   finite individual supports we can give simple predicate structures  that are closed  with respect to second-order definability.  The predicate structure containing, for every $n\geq 1$, exactly  those  predicates $ \alpha  \subseteq I^n$  with a finite individual support is a Henkin structure (cf.\,\,Example \ref{HenkinBasicMa}). This will be a direct consequence  of  Theorem \ref{HenkinStrIdeal}  because of certain relationships described by  properties such as  \ref{Eigenschaft3}.   
Let   $\G$ be a subgroup of $\G_1^I$  and $\pi\in \G$. For any subset $P\subseteq I$, let $ \pi P=\{\pi(\xi)\mid \xi\in P\}$  and, for any subgroup $\H$ of $\G$, let $\pi \H\pi^{-1}$ be  the subgroup of $\G$ given by $\pi \H\pi^{-1}=\pi\circ \H\circ \pi^{-1}=_{\rm df}\{\pi\circ \pi_0\circ \pi^{-1}\mid \pi_0\in \H\}$. 
\begin{enumerate}[label=\rm (P\arabic*),resume]
 \labelwidth0.7cm \leftmargin0cm \itemsep2pt plus1pt
\topsep1pt plus1pt minus1pt
\labelsep4pt \parsep0.5pt plus0.1pt minus0.1pt

\item\label{Eigenschaft3}   ${\G}( P)\subseteq {\rm sym}_{\G}(\alpha)$ implies
 ${\G}( \pi P)\subseteq \pi \G(P)\pi^{-1} \subseteq {\rm sym}_{\G}(\alpha^\pi)$.
\end{enumerate}

 \vspace{0.1cm}
 Thus, it would be  possible to  define Henkin structures with the help of normal filters as considered  in \cite{Jech03}  since each  normal filter is a  filter that is closed with respect to  all  transformations of  a subgroup $\H$ into a conjugate group $\pi\H\pi^{-1}$ by a permutation $\pi\in \G$.   In the following, we want to show 
that all predicates $\alpha$  over  $I$  for which ${\rm sym}_{\G}(\alpha)$    belongs to  such a normal filter    form a Henkin structure. Let us start with the definition of normal filters.

\section{Normal filter and permutation models} We follow \cite{1c}  where the fundamental  papers \cite{M38},  \cite{M39}, and   \cite{FA22a} are cited  in this context. First, we describe  a  construction method  for defining Henkin structures by using normal filters.  By \cite{Jech03}, Specker  \cite{Specker} extended the Fraenkel-Mostowski method by introducing the systems of subgroups of an automorphism group --- the so-called filters --- to construct models of first-order set theory.  By analogy with   \cite{Jech73} where Jech used normal filters  for constructing models of ZFA, we will use this kind of filters  $\F$ on  a group $\G$ of permutations of a non-empty individual  domain $I$ (for the used definition see \cite{TJJ}) for constructing  permutation models  of HPL that we will denote by $\S(I,\G,\F)$.
\begin{defi}[Normal filter]\label{DefNormFilt}Let $I$ be any non-empty individual  domain and $\G$ be a subgroup of $\G_{\sf 1} ^I$.  A system  $\F$ of subgroups of $\G$ that has the following properties (where $\H,\H_1,\H_2$ are subgroups of $\G$) is called a {\em normal filter on }  $\G$. 
\begin{enumerate}[label=(\roman*)]
 \labelwidth0.7cm \leftmargin0cm \itemsep0pt plus1pt
\topsep 1pt plus1pt minus1pt
\labelsep4pt \parsep0.5pt plus0.1pt minus0.1pt 
\item $ \G \in \F$;
\item if $\H_{1} \in \F$ and $\H_{2}$ is a subgroup of $ \G$ with  $\H_{2} \supseteq \H_{1}$,
      then  $\H_{2}$ is  in $\F$;
\item if $\H_{1},\H_{2} \in \F$, then the intersection $\H_{1} \cap \H_{2}$ is in $ \F$;
\item if $\H \in \F$ and $\pi \in {\G}$, then the conjugate subgroup 
      $\pi \H \pi^{-1}$ is in $ \F$;
\item if $P$ is a finite subset of $I$, then  $\G(P)$ is in  $\F$.
\end{enumerate}\end{defi}

\begin{remark} The  known definitions of such systems of subgroups  --- used  for constructing models of first-order set theory ---  differ from author to author (cf. \cite{Jech73,Howard,Felgner,Specker, TJJ}).
\end{remark}

\begin{example}[The normal filters $\F_{\sf \! 0}(I,\G)\!$ and $\F_{\sf \! 1}(I,\G)$] \hfill For each non-\linebreak empty individual  domain $I$ and any group $\G\subseteq \G_{\sf 1} ^I$, the  following systems are normal filter on $\G$.
\[\F_{\sf \! 1}(I, \G)=_{\rm df}\{\H \subseteq \G\mid \H \mbox{ \rm  is a subgroup of }\G\}\]
\[\F_{\sf \! 0}(I, \G)=_{\rm df}\{\H \in \F_{\sf \! 1}(I, \G)\mid (\exists P\subseteq I)( \mbox{\rm P is a finite set} \,\,\,\&\,\, \H\supseteq \G(P))\}\]
\end{example}
Asser  \cite{1c} considered  the simplest normal filters  $\F_{\sf \! 0}(I, \G)$ that can be generated by  all subgroups $\G(P)$   for which $P$  is a finite subset of $I$ in order to motivate the introduction of the notion {\em normal filter} for constructing Henkin structures.

\begin{defi}[Symmetric predicates] 
 Let $I$ be a non-empty set of individuals, let $\G$ be any   subgroup of  $\G_{\sf 1} ^I$,  and let $\F$ be any normal filter  on $\G$. Let $\alpha $ be a predicate in  ${\rm pred}(I)$. Then,  we say that $\alpha$  is {\em symmetric with respect to $\F$}  if the  symmetry subgroup ${\rm sym}_{\G}(\alpha)$ is in $ \F$.
\end{defi} 

\begin{defi}[The$\!$ second-order permutation model  $\S(I,\G,\F)$]\label{DefiPermModelByFilter}\hfill
 Let \linebreak $I$ be a  non-empty set,    $\G$ be any  subgroup of  $\G_{\sf 1} ^I$, and $\F$  be a normal filter  on $\G$. Then, for  any  $n\geq1$, let $J_{n}(I,\G,\F)=\{\alpha \in {\rm pred}_n (I)\mid {\rm sym}_{\G}(\alpha)\in \F\}$ and let  $\S(I,\G,\F)$ be the predicate structure $ (J_{n})_{n\geq 0}$ with the domains $J_{0}=I$ and $J_{n}=J_{n}(I,\G,\F)$.
\end{defi} 

Thus, for every $n\geq 1$,  $J_{n}(I,{\G,\F})$ is  the set of all $n$-ary predicates  that are symmetric with respect to $\F$.  The following closeness property   shows the usefulness of the normal filters (cf.\,Asser  \cite[Hilfssatz 3, p.\,113]{1c}).

\begin{proposition}\label{Closen_Filt} For any non-empty set $I$, any subgroup $\G$ of  $\G_{\sf 1} ^I$, and every normal filter $\F$ on $\G$,  the structure $\Sigma(I,{\G,\F})$ is closed with respect to all permutations $\pi \in \G$.
\end{proposition}

 \vspace{0.1cm}\noindent {\bf Proof.}  For any non-empty set $I$, any subgroup $\G$ of  $\G_{\sf 1} ^I$, and every normal filter $\F$ on $\G$, and  any $\alpha \in J_n(I,{\G,\F})$, there is an $\H\in \F$ with $\H\subseteq {\rm sym}_{\G}(\alpha)$. Let $\pi\in \G$. Analogously to  the part  (b) of the proof of property  \ref{Eigenschaft3},  we get  $ \pi \H\pi^{-1} \subseteq {\rm sym}_{\G}(\alpha^\pi)$ and thus $ {\rm sym}_{\G}(\alpha^\pi)\in \F$ and consequently $\alpha^\pi \in J_n(I,{\G,\F})$ by definition.
\qed

 \vspace{0.1cm} \noindent Now, we will characterize the stabilizers for the truth values of formulas by the following generalization of  Hilfssatz 4 in {\normalfont \cite[ p.\,113]{1c}}. In particular, we can say that  the subgroup   $ {\rm sym}^*_{\G}(\alpha)$  is  {\em the  stabilizer subgroup of  $\G^*$ with respect to a formula $H$ under an  assignment  $f$ in $\Sigma(I,{\G,\F})$} if $\alpha=\alpha_{\S,H,\mbmss{x}, f}$.   Proposition \ref{Closen_Filt}, together with Lemma \ref{HSAsser}\,(\ref{HSAsser2}), implies  the following Proposition \ref{AsserHS4} that can be proved by using the observations discussed for finite supports in  \cite{1c}. 

\setcounter{equation}{0}

Let us prove the following proposition used in \cite{Gass24A}.

\begin{proposition}[Stabilizers for formulas]\label{AsserHS4} Let  $I$ be a non-empty set, $\G$ be a subgroup of  $\G_{\sf 1} ^I$,  and $\F$  be a normal filter on $\G$. Let    $\S=\Sigma(I,{\G,\F})$ and $f\in {\rm assgn}(\S)$. Moreover,   let $H(\mbm{x})$ be in ${\cal L}^{(2)}_{\mbmss{x}}$    and $r,s\geq 0$. If $s>0$, then let $x_{i_1},\ldots,x_{i_s}$  ($ i_1,\ldots,i_s>n$) be  the individual  variables  that  may  additionally occur     free in $H(\mbm{x})$. If $r>0$,  then  let  $A_{j_1}^{n_1}, \ldots, A_{j_r}^{n_r}$  be the predicate variables that may occur     free in $H(\mbm{x})$.    Let $H(\mbm{x})$ does not contain further  free variables. Let $\H_{r+1}= \G(\{f(x_{i_1}) , \ldots, f(x_{i_s})\})$ if $s>0$ and $\H_{r+1}= \G$ if $s=0$. If $r>0$, then, for every  $k\in\{1,\ldots,r\}$, let        $A_k$ stand for  $A_{j_k}^{n_k}$ and let $\H_k $ be a subgroup in $\F$ such that  ${\rm sym}_{\G}(\alpha_k)\supseteq \H_k$  holds for $\alpha_k=f(A_k)$.    For each $l\in \{1,\ldots,r+1\}$,   let $\G_l$ be  the subgroup   in   $\F$ given  by  $\G_l= \H_{l}\cap \cdots \cap\H_{r+1}$.  Then, we have   
\begin{equation}\label{Stab_1}\S_f(H(\mbm{x},A_1,\ldots, A_{m}))=\S_{f\langle {\mbmty{x}\atop (\pi(f(\mbmty{x}))} {A_1\atop  (\alpha_1)^ \pi} {\cdots\atop\cdots} {A_{m}\atop (\alpha_{m})^ \pi}  \rangle}(H(\mbm{x},A_1,\ldots, A_{m}))\end{equation}
for each each $m\in \{1,\ldots, r\}$  and $\pi \in \G_{m+1}$ if   $r>0$ and
 \begin{equation}\label{Stab_2}{\rm sym}_{\G}(\alpha_{\S,H,\mbmss{x}, f})\supseteq \G_1\end{equation} in any case where $r\geq 0$.
\end{proposition}

\noindent{\bf Proof.} Let $\pi\in \G_1$. Then,  $\pi\in \G$ and $\pi^{-1}\in \G$. Thus,  we have $\alpha_{\S,H,\mbmss{x}, f^\pi}(\mbm{\xi})=\alpha^\pi_{\S,H,\mbmss{x}, f}(\mbm{\xi})$  for all $\mbm{\xi}\in I^n$ by Lemma \ref{HSAsser} (\ref{HSAsser2}). Clearly, for $r\geq 1$ and $l\in\{1,\ldots,r\}$, $\alpha_l\in  J_{n_l}(I,{\G,\F})$ implies $(\alpha_l)^ \pi\in  J_{n_l}(I,{\G,\F})$ by Proposition \ref{Closen_Filt}. Therefore,  for any $\mbm{\xi}\in I^n$,  by definition and Proposition \ref{MeaningOfFreeVar}, there holds 

\[\begin{array}{llr}\!\!\alpha_{\S,H,\mbmss{x}, f}(\mbm{\xi})\!\!&=\S_{f\langle {\mbmty{x}\atop\mbmty{\xi}}\rangle}(H(\mbm{x}))&(\mbox{by definition})\\&= \S_{f \langle {x_{i_1}\atop f(x_{i_1})} {\cdots \atop \cdots} {x_{i_s}\atop f(x_{i_s})}  {A_1\atop \alpha_1}{\cdots \atop \cdots}  {A_r\atop \alpha_r}\rangle\langle {\mbmty{x}\atop\mbmty{\xi}}\rangle}(H(\mbm{x}))\\

&=\S_{f \langle {x_{i_1}\atop f^{\pi}(x_{i_1})} {\cdots \atop \cdots} {x_{i_s}\atop f^{\pi}(x_{i_s})}  {A_1\atop (\alpha_1)^{\pi}}{\cdots \atop \cdots}  {A_r\atop( \alpha_r)^{\pi}}\rangle\langle  {\mbmty{x}\atop\mbmty{\xi}}\rangle}(H(\mbm{x}))&(\mbox{since } \pi\in \G_1)

\\&=\S_{f ^{\pi}\langle {x_{i_1}\atop f^{\pi}(x_{i_1})} {\cdots \atop \cdots} {x_{i_s}\atop f^{\pi}(x_{i_s})}  {A_1\atop (\alpha_1)^{\pi}}{\cdots \atop \cdots}  {A_r\atop (\alpha_r)^{\pi}}\rangle\langle  {\mbmty{x}\atop\mbmty{\xi}}\rangle}(H(\mbm{x}))\!\!\!\!\!&\\&&\!\!\!\!\!\!\!\!\!\!\!\!\!\!\!\!\!\!\!\!(\mbox{by Proposition }\ref{MeaningOfFreeVar})

\\&=\S_{(f \langle {x_{i_1}\atop f(x_{i_1})} {\cdots \atop \cdots} {x_{i_s}\atop f(x_{i_s})}  {A_1\atop \alpha_1}{\cdots \atop \cdots}  {A_r\atop \alpha_r} \rangle )^{\pi} \langle {\mbmty{x}\atop\mbmty{\xi}}\rangle}(H(\mbm{x}))

\\&=\S_{f  ^{\pi} \langle {\mbmty{x}\atop\mbmty{\xi}}\rangle}(H(\mbm{x}))\\
&=\alpha_{\S,H,\mbmss{x}, f^\pi}(\mbm{\xi})&(\mbox{by definition})

\\&=\alpha^\pi_{\S,H,\mbmss{x}, f}(\mbm{\xi})&\!\!\!\!\!\!(\mbox{by Lemma }\ref{HSAsser} (\ref{HSAsser2})).
\end{array}
\]
Consequently, there hold $\alpha^\pi_{\S,H,\mbmss{x}, f}=\alpha_{\S,H,\mbmss{x}, f}$ for any $\pi\in \G_1$ and thus the inclusion (\ref{Stab_2}). By Lemma \ref{HSAsser} (\ref{HSAsser1}), we  can use $\S_{f}(H)=\S_{f^\pi}(H)$ for $H=H(\mbm{x},A_1,\ldots,A_m)$. For  $\pi \in \G_{m+1}$, we get 

\[\begin{array}{llr}\S_f(H)\\

\hspace*{1cm}= \S_{f \langle{\mbmty{x}\atop f(\mbmty{x})} {x_{i_1}\atop f(x_{i_1})} {\cdots \atop \cdots} {x_{i_s}\atop f(x_{i_s})}  {A_1\atop \alpha_1}{\cdots \atop \cdots}  {A_m\atop \alpha_m} {A_{m+1}\atop \alpha_{m+1}}{\cdots \atop \cdots}  {A_r\atop \alpha_r} \rangle}(H)\\

\hspace*{1cm}= \S_{(f \langle{\mbmty{x}\atop f(\mbmty{x})} {x_{i_1}\atop f(x_{i_1})} {\cdots \atop \cdots} {x_{i_s}\atop f(x_{i_s})}  {A_1\atop \alpha_1}{\cdots \atop \cdots}  {A_m\atop \alpha_m} {A_{m+1}\atop \alpha_{m+1}}{\cdots \atop \cdots}  {A_r\atop \alpha_r} \rangle)^\pi}(H)\\

\hspace*{1cm}= \S_{f( \langle{\mbmty{x}\atop f(\mbmty{x})} {x_{i_1}\atop f(x_{i_1})} {\cdots \atop \cdots} {x_{i_s}\atop f(x_{i_s})}  {A_1\atop \alpha_1}{\cdots \atop \cdots}  {A_m\atop \alpha_m} {A_{m+1}\atop \alpha_{m+1}}{\cdots \atop \cdots}  {A_r\atop \alpha_r} \rangle)^\pi}(H)\\

\hspace*{1cm}= \S_{f \langle{\mbmty{x}\atop \pi(f(\mbmty{x}))} {x_{i_1}\atop \pi(f(x_{i_1}))} {\cdots \atop \cdots} {x_{i_s}\atop \pi( f(x_{i_s}))}  {A_1\atop (\alpha_1)^\pi}{\cdots \atop \cdots}  {A_m\atop (\alpha_m)^\pi} {A_{m+1}\atop (\alpha_{m+1})^\pi}{\cdots \atop \cdots}  {A_r\atop (\alpha_r)^\pi} \rangle}(H)\\

\hspace*{1cm}= \S_{f \langle{\mbmty{x}\atop \pi(f(\mbmty{x}))} {x_{i_1}\atop f(x_{i_1})} {\cdots \atop \cdots} {x_{i_s}\atop f(x_{i_s})}  {A_1\atop (\alpha_1)^\pi}{\cdots \atop \cdots}  {A_m\atop (\alpha_m)^\pi} {A_{m+1}\atop \alpha_{m+1}}{\cdots \atop \cdots}  {A_r\atop \alpha_r} \rangle}(H)\\

\hspace*{1cm}= \S_{f \langle{\mbmty{x}\atop \pi(f(\mbmty{x}))}  {A_1\atop (\alpha_1)^\pi}{\cdots \atop \cdots}  {A_m\atop (\alpha_m)^\pi}  \rangle}(H)&\hspace{0.85cm}\Box
\end{array}
\] 

\noindent Equation (\ref{Stab_2}) in Proposition \ref{AsserHS4}  means that, for  each predicate  definable  in a predicate structure $\Sigma(I,{\G,\F})$, the  symmetry subgroup of $\G$   with respect  to this predicate   belongs to the considered filter $\F$. This  guarantees the closeness of $\Sigma(I,{\G,\F})$ with respect to the  definability  by means of second-order formulas. This means that  $\Sigma(I,{\G,\F})$ is closed under  $Att_\S^n$ for all $n\geq 1$. A consequence is   that  the  structures constructed are models of the axioms  given in (A4).

\begin{theorem}[Asser {\normalfont \cite[p.\,114]{1c}}]\label{AsserHenkin}  For any non-empty set $I$ of individuals, any subgroup $\G$ of  $\G_{\sf 1} ^I$, and each normal filter $\F$ on $\G$, the  predicate structure $\S(I,\G,\F)$ is a Henkin structure.
\end{theorem}

Specker  used  filters  on subgroups of groups of automorphisms of a set so that the relation $\in$ stays invariant with respect to these  automorphisms.  Mostowski  (cf.\,\,also \cite{Jech73})   used order-preserving permutations (cf.\,\,\cite{M39}) of rational numbers    represented by sets $\Lambda_n$ (defined in  \cite{M39}) for showing the independence of the well-ordering theorem from the principle of the total ordering  \cite{M39}. Thus, we can say he applied a group defined by  automorphisms of $(\mathbb{Q};\emptyset;\emptyset;<)$ whose extensions preserve  the relation $\in$.  Here,     we can use  filters on various  groups of automorphisms of suitable mathematical  structures. In general, we will take  suitable predicate basic structures  which can  be derived from structures  containing   several mathematical objects such as functions, relations, systems of sets or  only individuals.  Every predicate basic structure $\S({\sf S})$  derived from a first-order basic structure ${\sf S}$ in $ {\sf struc}^{(1)}(I)$ contains all predicates $\alpha$  for which  the corresponding relation $\widetilde \alpha$ belongs to  the underlying structure ${\sf S}$   or   the  graph  of a function $\mbm{g}:I^{p}\to I$ in ${\sf S}$  given by ${\rm graph}(\mbm{g})=_{\rm df}\{( \mbm{\xi}\,.\,\eta)\mid \mbm{\xi}\in I^{p} \,\,\&\,\, \mbm{g}( \mbm{\xi})= \eta\}$    is the corresponding relation $\widetilde \alpha\subseteq  I^{p+1}$ where  $(\mbm{\xi}\,.\,\eta)$ is here   the tuple $(\xi_1, \ldots, \xi_{p}, \eta)$.

\section{Normal ideals and permutation models} In analogy to the construction of models  of  ZFA considered in  \cite[p.\,47]{Jech73}, a lot of 
normal filters  on a  group $\G$ of permutations of a domain $I$ can be defined by  using any normal ideal $\I\subseteq {\cal P}(I)$.   ${\cal P}(I)$ is here the power set that  contains all subsets of $I$ within the framework of ZF.   Its  existence is guaranteed by the axiom of the power set.  In \cite{Gass94}, we extended the notion of normal ideal and  introduced systems  $\I\subseteq {\cal P}(I\cup {\rm pred}(I))$ that are  closed  with  respect to the left group action of $\G^*$.  For any $P\in {\cal I}$,   $(\G(P))^*$  is a subgroup of $\G^*$ that is  included in each stabilizer  subgroup of $\G^*$  with respect to  one  of the  individuals $\xi$ in $P$ and the  relations $\alpha$ in $P$, respectively. For all $\pi \in {\G}$ and  for all $P\subseteq I\cup {\rm pred}(I)$, let
 \[\pi P=_{\rm df}\{\pi(\xi)\mid \xi\in P \,\,\& \,\, \xi \in I \}      \cup \{\alpha^\pi\mid  \alpha \in P \,\,\& \,\, \alpha \in {\rm pred}(I)\}.\] 
This means  that we have  inclusions  such as in \ref{Eigenschaft3}.
\begin{enumerate}[label=\rm (P\arabic*),resume]
 \labelwidth0.7cm \leftmargin0cm \itemsep2pt plus1pt
\topsep1pt plus1pt minus1pt
\labelsep4pt \parsep0.5pt plus0.1pt minus0.1pt

\item\label{Eigenschaft4}   ${\G}( P)\subseteq {\rm sym}_{\G}(\alpha)$ implies
 ${\G}( \pi P)\subseteq \pi \G(P)\pi^{-1} \subseteq {\rm sym}_{\G}(\alpha^\pi)$.
\end{enumerate}
We can prove them analogously to the inclusions in  \ref{Eigenschaft3}. 
Thus,  the extended ideals  of second order are suitable for the definition of further normal filters. These properties motivate us to give Definition \ref{DefExtNormIdeal}.

 \vspace{0.1cm}

 \begin{defi}[Extended normal ideal]\label{DefExtNormIdeal}  For any domain $I$ and any subgroup $\G$ of  $\G_{\sf 1} ^I$,  a system $\I$  of subsets of $I\cup {\rm pred}(I)$ is  an {\em extended normal ideal on $ I$ with respect to $\G$}
 if $\I$ has the following five properties. 

 \vspace{0.1cm}

\begin{enumerate}[label=(\roman*)]
 \labelwidth0.7cm \leftmargin0cm \itemsep0pt plus1pt
\topsep1pt plus1pt minus1pt
\labelsep4pt \parsep0.5pt plus0.1pt minus0.1pt 
\item $ \emptyset \in \I$;
\item if $P_{1} \in \I$ and $P_{2} \subseteq P_{1}$,  then $P_{2} \in \I$;
\item if $P_{1},P_{2} \in \I$, then $P_{1} \cup P_{2}\in \I$;
\item\label{Propertyiv} if $P \in \I$ and $\pi \in {\G}$, then    $\pi P \in \I$;
\item if $P$ is a finite subset of $I$, then $P \in \I$.
\end{enumerate}
 \vspace{0.1cm}

If $\I$  contains only subsets of $I$, then it  is  a  ({\em simple}) {\em  normal ideal on $ I$ with respect to $\G$}. In general we omit the words {\em simple} and  {\em extended} and call  the extended normal ideals also   {\em normal ideals} ({\em of second order}).  
\end{defi}

The property {\it\ref{Propertyiv}}  means that   $\I$ is closed under the left action of $\G^*$.

\begin{example}[The simple normal ideals $\I_{\sf 0}^I$ and $\I_{\sf 1}^I$]\label{ExamNormId1a}\hfill  Let $I$  be any non-\linebreak empty set. Then, we can consider the following systems. 
\[\I_{\sf 0}^I=\{P\subseteq I\mid P \mbox{ is a finite set}\}\]
\[\I_{\sf 1}^I={\cal P}(I)\]

 \noindent  With respect to each group $\G\subseteq \G_{\sf 1} ^I$, $\I_{\sf 0}^I$ is the simplest normal ideal on $I$ and  the simple normal ideal $\I_{\sf 1}^I$ is the maximal extended ideal $\subseteq {\cal P}(I)$. 
\end{example}
\begin{example}[The extended normal ideal $\I_{\sf 0}^{I\cup P_0}$]\label{ExamNormId1b} Let $I$  be a non-empty set, $\G$ be a subgroup of  $\G_{\sf 1} ^I$,  and $P_0$ be a non-empty subset of  ${\rm pred}(I)$  that is closed with respect to all $\pi\in \G$, which means that  $\alpha^\pi $ belongs to $P_0$  for each predicate  $\alpha \in P_0$ and each $\pi\in \G$.
Moreover, let

$\I_{\sf 0}^{I\cup P_0}=\{P\subseteq I\cup P_0 \mid P \mbox{ is a finite set}\}$.

 \noindent  Then, $\I_{\sf 0}^{I\cup P_0}$ is an   extended ideal $\subseteq {\cal P}(I\cup P_0)$ with respect to $\G$. 
\end{example}

We will use  the possibility to define  normal filters  by means of  normal ideals.

 \vspace{0.1cm}

\begin{example}[The normal filter $\F_{\sf \! 0}(I,\G,\I)$ induced by  ideal $\I$]\label{ExamNormId2}  \hfill For \linebreak any subgroup $\G$ of  $\G_{\sf 1} ^I$ and any extended normal ideal $\I$ on $ I$ with respect to  $\G$, the  following system is a normal filter on $\G$  induced by  the  ideal $\I$. 
\[\F_{\sf \! 0}(I, \G,\I)=_{\rm df}\{\H \in \F_{\sf \! 1}(I, \G)\mid (\exists P\in \I) (\H\supseteq \G(P))\}\] 
Thus, there holds
\[\F_{\sf \! 0}(I, \G,\I_{\sf 0}^I)= \F_{\sf \! 0}(I, \G).\] 
By property \ref{Eigenschaft4}, the closeness of the normal ideals $\I$ with respect to applying a permutation $\pi\in \G$ on a set $P\in \I$ as described by {\it\ref{Propertyiv}} in Definition \ref{DefExtNormIdeal}  implies that the induced filters $\F$ are closed with respect to  the  transformation of  any subgroup $\H\in \F$ into a conjugate group $\pi\H\pi^{-1}$ by  $\pi$.
\end{example}

This implies the possibility to define Henkin structures by extended normal ideals directly. 

\begin{defi}[The  permutation model   $\!\S(I,\G,\I)$]$\!$Let $I$ be
a non-empty set. Let     $\G$ be any  subgroup of  $\G_{\sf 1} ^I$ and $\I$  be an extended normal ideal on  $I$  with respect to  $\G$. Then, for  any  $n\geq1$, let  
\[J_{n}(I,\G,\I)=\{\alpha \in {\rm pred}_n(I)\mid ( \exists P \in {\I})( {\rm sym}_{\G}(\alpha)\supseteq{\G}( P))\}\]
and let $\S(I,\G,\I)$ be the predicate structure $ (J_{n})_{n\geq 0}$ with the domains $J_{0}=I$ and $J_{n}=J_{n}(I,\G,\I)$. 
\end{defi}

\begin{theorem}[Ga\ss ner {\normalfont \cite[p.\,542]{Gass94}}]\label{HenkinStrIdeal} For every non-empty individual domain $I$, any group  $\G$ of permutations of $I$, and each (extended) normal ideal $\I$ on $I$ with respect to  $\G$,  the  structure $\S(I,\G,\I)$  is a Henkin-Asser structure.
\end{theorem}  \vspace{0.1cm}

\begin{example}[The   basic Henkin structure  $\S_{\sf 0}^I$ of second order]\label{HenkinBasicMa} \hfill For \linebreak  any non-empty set  $I$,   the predicate structure  $\S(I, \G_{\sf 1} ^I,\I_{\sf 0}^I)$  is the Henkin structure   that we denote by $\S_{\sf 0}^I$.
 $\S_{\sf 0}^I$ is called the  {\em basic  Henkin structure over the domain $I$} of second order.
\end{example}

The basic Fraenkel model  of   second order, $\S_0^\bbbn$   can be  defined  in analogy with the basic Fraenkel model (cf.\,\,\cite[p.\,48]{Jech73}) (see also \,\,\cite[Section V.1, pp.\,47--49]{Gass84}, \cite[Section 7, pp.\,544--545]{Gass94}). 

\begin{example}[The standard structure  $\S_{\sf 1}^I$ of second order]\label{HenkinBasicMb} \hfill For  any \linebreak non-empty set  $I$,  let $\G_{\sf 0} ^I$ contain only the identity permutation of $I$. Then,  $\G_{\sf 0} ^I$ is a subgroup of $\G_{\sf 1} ^I$, and the Henkin structure $\S(I, \G_{\sf 0} ^I,\I_{\sf 0}^I)$  is the standard structure  $\S_{\sf 1}^I$.
\end{example} 

In \cite{Gass84} (and so in \cite{Gass94}), various Henkin  structures are constructed by using suitable groups and different normal ideals in order to compare several variants of the  Axiom of Choice and some other statements.  Asser introduced the ordered Mostowski model of second order in analogy to the ordered Mostowski model (cf.\,\,\cite[p.\,50]{Jech73}) as follows. The construction can be done in ZFC.

 \vspace{0.1cm}

\begin{example}[The Mostowski-Asser model  $\S_1$ {\normalfont \cite[p.\,114]{1c}}]\label{ModelMost}  \hfill Here, let $\bbbn$ \linebreak be as usually the least non-zero limit ordinal $\omega$ (for more details, see \cite[p. 20]{Jech03}). Let  $\bbbq$ be the set  of all rational numbers ordered by  the usual ordering $\leq$.   The numbers in $\bbbq$ are in general defined by  the classes represented  by $0$ and $\frac{n}{m}$ and $-\frac{n}{m}$, respectively,   where  $n$ and $ m$ are coprime  numbers in $\omega$  with  $n\not=0$ and $m\not=0$. They can also be represented by sets  such as  $\Lambda_n$ in  \cite{M39}. The restriction of $\leq$ on $\mathbb{Q}$ is the linear ordering  $\mbm{r}$ given by $\mbm{r}=\{(q_1,q_2)\in \bbbq\times\bbbq\mid q_1\leq q_2\}$ which, for the sake of simplicity, we also denote by $\leq$.  We call the Henkin structure $\S(\bbbq,{\rm auto}(\bbbq;\emptyset; \emptyset; \leq),\I_{\sf 0}^{\bbbq})$ the {\em [ordered]  Mostowski-Asser  model [of second order]} and denote it by $\S_1$.  Note that $\S_1$ could also be defined  by $\S_1=\S(\bbbq,{\rm auto}(\bbbq;\emptyset; \emptyset; \leq),\I_{\sf 0}^{\bbbq\cup \{\leq\}})$ and so on. 
\end{example}

 \vspace{0.1cm}

The following model is derived from  a permutation model of  ZFA  that Jech  called the second Fraenkel model (cf.\,\,\cite[p.\,48]{Jech73}).  It was introduced by Fraenkel  (cf.\,\cite{FA22a,FA37}) by using also  an idea of  Zermelo. An extension of the following  model is introduced in \cite[Section V.8, pp.\,92--96]{Gass84}.

\begin{example}[The   Fraenkel-Zermelo model  of second order]\label{Model_sec_Frae} \hfill
 For \linebreak each $n\in \bbbn$,  let $\mbm{t}_n = \{((n,0),(n,1)),((n,1),(n,0))\}$. Then,  each $\mbm{t}_n $ is a  binary relation on $I=\bbbn\times\{0,1\}$.   The structure $\S(I,{\rm auto}(I;\emptyset;\emptyset;(\mbm{t}_n)_{n\geq 0}),\I_{\sf 0}^I)$ is a Henkin structure.
\end{example}

\section{The use of permutation models of second order}

In \cite{Gass84}, 14 Henkin-Asser structures M(1), \ldots, M(14) are investigated.   

The first 12 of these structures considered also in the following overviews are permutation models.  The propositions are called {\sf Satz} as in \cite{Gass84}  and the numbers given in the last column of the tables refer to the propositions in \cite{Gass84}. ${\rm M(1)}$ is  here the basic Fraenkel model  $\S_0$ of second order. 
For the definitions, see \cite{Gass84} and \cite{Gass94}. For further overviews see \cite[pp. 31--32]{Gass84} and \cite{Gass94}.

 \vspace{0.1cm}

Note, that Stephen Mackereth considered  some similar problems  in \cite{Mackereth}.

 \vspace{0.3cm}

Let $2\S_0=\S(I, \G,\I_0^I)$ with $I=\{0,1\}\times \bbbn$ and  
   \[\G\!=\!\{\pi\in \G_1^I\mid  \pi^ *|_{\S_0^{\{0\}\times \bbbn}}\in {\rm auto}^*(\S_0^{\{0\}\times \bbbn}) \,\,\&\,\pi^ *|_{\S_0^{\{1\}\times \bbbn}}\in {\rm auto}^*(\S_0^{\{1\}\times \bbbn})  \}\] which means that $\pi \in \G$ implies 
  \[\xi \!\in \{0\}\!\times \bbbn\Rightarrow  \pi(\xi)\!\in \{0\}\!\times \bbbn\mbox{ \,\, and \,\, } \xi \!\in \{1\}\!\times \bbbn\Rightarrow \pi(\xi)\!\in \{1\}\!\times \bbbn\]
 for all $\xi \in I$. For ${\cal G}\subseteq {\cal L}^{(2)}$, $\S\models {\cal G}$ means that  $\S\models H$ holds for all $H\in {\cal G}$.

\begin{remark}[A possible  answer for a problem  in {\cite[p.\,\,16]{Mackereth}}]\hfill
By \cite{Gass24A} \linebreak and \cite{Gass24B}, $\S_0\models {\rm HAC}$ holds.  The same holds also for $2\S_0$.
We have $2\S_0\models {\rm HAC}$. Consequently,  $(^hax^{(2)}\cup {\rm HAC})\vdash TR^1$ and $(^hax^{(2)}\cup {\rm HAC})\vdash LO^1$ do not hold  by Proposition  \ref{PROP_TR_LO}.  

 \end{remark}

  \begin{proposition}\label{PROP_TR_LO}\hfill 
  
  \begin{tabular}{lll}
  
 {\rm (1)} &$2\S_0\models \neg TR^n$ & for all $n\geq 1$.\\
  
 {\rm (2)} & $2\S_0\models \neg LO^n$ & for all $n\geq 1$.\\
  
 {\rm (3)} & $2\S_0\models choice_h^{n,m}$ & for all $n,m\geq 1$.\\
  \end{tabular}
  
  \end{proposition}
  
\noindent  {\bf Proof.} Property (3)  of  Proposition \ref{PROP_TR_LO}   follows from Proposition 4.14 in \cite{Gass24B}.   Property (2) follows from \cite[Section V, Satz 1.5]{Gass84}.   Property (1) holds because  no  injective mapping   of an  infinite subset of $\{0\}\times \bbbn$ into an infinite subset of $\{1\}\times \bbbn$ belongs to $\widetilde{J_{2n}}$  for $J_{2n}=J_{2n}(I, \G,\I_0^I)$ (for $n=1$). \qed
  
  More general, if ZFC is consistent, then we have the following.
  
  \begin{theorem} For all $n\geq 1$, $TR^n$ and $LO^n$ are  {\rm HPL}-independent of $choice_h^{(2)}$.
  \end{theorem}

\newpage

\begin{overview}{\sf  The models from $M(1)$ to $M(7)$ (a summary)}
\nopagebreak 

\noindent
{
\small
\begin{tabular}{|cc|ll|c|}
\hline
Section$\!\!\!$&Model &Positive results&Negative results&Satz\\ 
 \hline 

\cite[Sect. V.1]{Gass84}&$M(1)$&$M(1)\models TR^1$&& 1.2\\  

 &&$M(1)\models AC^{1,1}$&& 1.3\\  

 &&$M(1)\models LW^1$&& 1.4\\  

  &&&$M(1)\models \neg LO^1$ & 1.5\\    \hline

\cite[Sect. V.2]{Gass84}&$M(2)$&$M(2)\models AC^{1,1}$&& 2.5\\  

 &&$M(2)\models LO^1$&& 2.6\\  

  &&&$M(2)\models \neg TR^1$& 2.7\\  

  &&&$M(2)\models \neg LW^1$& 2.8\\  

  &&&$M(2)\models \neg KW$-$AC^{2,1}$& 2.9\\  

  &&&$M(2)\models \neg MC(4)_{LO}^{2,1}$& 2.10\\     \hline

\cite[Sect. V.3]{Gass84}&$M(3)$&$M(3)\models TR^1$&& 3.5\\  

 &&$M(3)\models DC^1$&& 3.6\\  

 &&$M(3)\models LO^1$&& 3.7\\  

 &&&$M(3)\models \neg AC^{1,1}_*$& 3.8\\  

 &&&$M(3)\models \neg MC(4)^{1,1}_{LO}$& 3.9\\  

 &&&$M(3)\models \neg KW$-$AC^{1,1}$& 3.10\\  

 &&&$M(3)\models \neg LW^1$& 3.11\\   

 &&$M(3)\models AC^{1,1}_\omega$&& 3.12\\    \hline

\cite[Sect. V.4]{Gass84}&$M(4)$&$M(4)\models TR^1$&& 4.2\\  

 &&&$M(4)\models \neg TR^2$& 4.3\\    \hline

\cite[Sect. V.5]{Gass84}&$M(5)$&$M(5)\models ZL^1$&& 5.2\\  

 &&$M(5)\models LW^1$&& 5.3\\  

 &&&$M(5)\models \neg MC(4)^{1,1}_{\omega}$& 5.4\\  

 &&&$M(5)\models \neg MC(4)^{1,1}_{LO}$& 5.5\\  

 &&&$M(5)\models \neg AC^{1,1}_*$& 5.6\\  

 &&&$M(5)\models \neg KW$-$AC^{1,1}$& 5.7\\  

 &&&$M(5)\models \neg DC^1$& 5.8\\   \hline

\cite[Sect. V.6]{Gass84}&$M(6)$&$M(6)\models GCH(2)^1$&& 6.2\\  

 &&$M(6)\models LW^1$&& 6.3\\  

 &&&$M(6)\models \neg LO^1$ & 6.4\\  

 &&&$M(6)\models \neg TR^1$& 6.5\\  

 &&&$M(6)\models \neg AC^{1,1}_*$& 6.6\\  

 &&&$M(6)\models \neg KW$-$AC^{1,1}$& 6.7\\    \hline

\cite[Sect. V.7]{Gass84}&$M(7)$&$M(7)\models GCH(1)^1$ && 7.2\\  

 &&$M(7)\models LO^1$&& 7.3\\  

 &&$M(7)\models AC^{1,1}_*$&& 7.4\\  

 &&&$M(7)\models \neg ZL^1$& 7.5\\  
 
  &&&$M(7)\models \neg DC^1$& 7.6\\

 \hline  

\end{tabular}
}

\end{overview}

\newpage

\begin{overview}{\sf  The models from  $M(7)$ to $M(14)$  (a summary)}

\nopagebreak 
\noindent
{\small
\begin{tabular}{|cc|ll|c|}
\hline
Section$\!\!\!$&Model &Positive results&Negative results&Satz\\ 
 \hline

\cite[Sect. V.7]{Gass84}

&&$M(7)\models AC^{1,1}_\omega$&& 7.7\\  

 &&&$M(7)\models \neg LW^1$& 7.8\\  

 &&&$M(7)\models \neg MC(2)_{LO}^{1,1}$& 7.9\\  \hline  

\cite[Sect. V.8]{Gass84}&$M(8)$&$M(8)\models MC(2)^{1,1}$&& 8.2\\  

 &&$M(8)\models LW^1$&& 8.3\\  

 &&&$M(8)\models \neg AC_*^{1,1}$& 8.4\\  

 &&&$M(8)\models \neg AC^{1,1}_\omega$& 8.5\\  

 &&&$M(8)\models \neg DC^1$& 8.6\\  

 &&&$M(8)\models \neg KW$-$AC^{1,1}$& 8.7\\  

 &&&$M(8)\models \neg ZL^1$& 8.8\\      \hline

\cite[Sect. V.9]{Gass84}&$M(9)$&$M(9)\models KW$-$AC^{1,1}\!\!\!$&& 9.3\\  

 &&&$M(9)\models \neg AC^{1,1}_*$& 9.4\\  

 &&&$M(9)\models \neg  ZL^1$& 9.5\\  

 &&&$M(9)\models \neg AC^{1,1}_\omega$& 9.6\\  

 &&&$M(9)\models \neg DC^1$& 9.7\\     \hline

\cite[Sect. V.10]{Gass84}$\!\!$&$\!\!$$M(10)$&$M(10)\models AC^{1,1}_*$&& 10.2\\  

 &&&$M(10)\models \neg MC(1)^{1,1}_{LO}$& 10.3\\  

 &&&$M(10)\models \neg KW$-$AC^{1,1}$& 10.4\\    \hline

\cite[Sect. V.11]{Gass84}$\!\!$&$\!\!$$M(11)$&$M(11)\models LO^1$&& 11.2\\  

 &&&$M(11)\models \neg MC(4)^{1,1}_{\omega}$& 11.3\\    \hline

\cite[Sect. V.12]{Gass84}$\!\!$&$\!\!$$M(12)$&$M(12)\models ZL^1$&& 12.2\\  

 &&&$M(12)\models \neg AC^{1,2}_*$& 12.3\\  

 &&&$M(12)\models \neg  ZL_2$& 12.4\\  

 &&$M(12)\models AC^{n,1}_*$&& 12.5\\    \hline  

\cite[Sect. V.13]{Gass84}$\!\!$&$\!\!$$M(13)$&$M(13) \models WO^1$&& 13.2\\    

 &&&$M(13)\models \neg GCH(4)^1$ & 13.3 \\   \hline

\cite[Sect. V.14]{Gass84}$\!\!$&$\!\!$$M(14)$&$M(14)\models DC^1$&& 14.2 \\  

 &&&$M(14)\models \neg ZL^1$& 14.3 \\ \hline

\end{tabular}

}
\end{overview}

\end{document}